\documentclass[amsmath,11pt]{article}

\usepackage{graphicx}    % standard LaTeX graphics tool
\usepackage{amssymb}
\usepackage{amsfonts}
\usepackage{bm}
\usepackage{amsmath}

\newtheorem{theorem}{Theorem}[section]
\newtheorem{lemma}[theorem]{Lemma}

\def\real {{I \!\! R}}

\def\bigM{ {\mathcal {M}}}
\def\bigF{ {\mathcal {F}}}

\def\g{\hat g}

\newcommand{\rf}[1]{(\ref{#1})}

\def\g{\gamma}

\def\({\left(}
\def\){\right)}

\def\be#1\ee{\begin{equation}#1\end{equation}}
\def\ba#1#2\ea{\begin{array}#1#2\end{array}}
\def\bgr#1\egr{{\allowdisplaybreaks\begin{gather}#1\end{gather}}}
\def\bma#1\ema{{\allowdisplaybreaks\begin{align}#1\end{align}}}

\def\oplem#1{\begin{lemma}\, {\rm #1}\, \it }
\def\cllem{\end{lemma}\rm \par }
\def\opthm#1{\begin{theorem}\, {\rm #1}\, \it }
\def\clthm{\end{theorem}\rm \par }

\newcommand{\fer}[1]{(\ref{#1})}
\newcommand{\bq}{\begin{equation}}
\newcommand{\eq}{\end{equation}}
\def\bqa{\begin{eqnarray}}
\def\eqa{\end{eqnarray}}
\def\bd{\begin{displaymath}}
\def\ed{\end{displaymath}}

\begin{document}
\title{On a kinetic model for a simple market economy}

\author{Stephane Cordier\thanks{Department of Mathematics and Applications, Mathematical Physics (MAPMO)UMR 6628
University of Orl\'eans and CNRS, 45067 Orl\'eans, France. \texttt{
cordier@labomath.univ-orleans.fr} }\,\and Lorenzo Pareschi\thanks{Department of Mathematics,
University of Ferrara, Via Machiavelli 35 I-44100 Ferrara, Italy, and Center for Modelling
Computing and Statistics, c/o Department of Economy, Institutions and Territory, University of
Ferrara, Via del Gregorio 13, I-44100 Ferrara, Italy. \texttt{pareschi@dm.unife.it}}\,\and
Giuseppe Toscani\thanks{Department of Mathematics, University of Pavia, via Ferrata 1, 27100
Pavia, Italy. \texttt{toscani@dimat.unipv.it}} }

\maketitle

%%%%%%%%%%%%%%%%%%%%%%%%%%%%%%%%%%%%%%%%%%%%%%%%%%%%%%%%%%%%%%%%%%%%%

\begin{abstract}
In this paper, we consider a simple kinetic model of economy
involving both exchanges between agents and speculative trading.
We show that the kinetic model admits non trivial quasi-stationary
states with power law tails of Pareto type. In order to do this we
consider a suitable asymptotic limit of the model yielding a
Fokker-Planck equation for the distribution of wealth among
individuals. For this equation the stationary state can be easily
derived and shows a Pareto power law tail. Numerical results
confirm the previous analysis.
\end{abstract}

{\bf Keywords.} Econophysics, Boltzmann equation, wealth and
income distributions, Pareto distribution.

\tableofcontents
%%%%%%%%%%%%%%%%%%%%%%%%%%%%%%%%%%%%%%%%%%%%%%%%%%%%%%%%%%%%%%%%%%%%%%%%%%

\section{Introduction}

Microscopic models of simple market economies have been recently
introduced by several authors \cite{AC, AC1, YD, IKR}. The main
idea is that an economic system composed by a sufficiently large
number of agents can be described using the laws of statistical
mechanics as it happens in a physical system composed of many
interacting particles. The details of the economical interactions
between agents characterize their wealth distribution.

The study of wealth distributions has a long history going back to
more then one century ago with the Italian sociologist and
economist Vilfredo Pareto which studied the distribution of income
among people of different western countries and found an inverse
power law \cite{Pa}. More precisely if $f(w)$ is the probability
density function of agents with wealth $w$ we have
\[
F(w)=\int_w^{\infty} f(w_*)\,dw_* \sim w^{-\mu}.
\]
Pareto mistakenly believed that power laws apply to the whole
distribution with a universal exponent $\mu$ approximatively equal
to $1.5$. Later, Mandelbrot, proposed a weak Pareto law that
applies only to high incomes \cite{Ma}.

Income data from developed countries show Pareto-like behaviors
with order $1$ exponents for large values of incomes (USA $\sim
1.6$, Japan $\sim 1.8-2.2$ \cite{YD}). It is common, in fact, that
$90\%$ of the total wealth is owned by only $5\%$ of the
population. According to data, across the full range of income, we
should expect the probability density function $f(w)$ to increase
at low income, reach a maximum and, finally decrease with
increasing wealth \cite{YD}.

In our discussion we use the terms distribution of income (money)
and distribution of wealth interchangeably (although the two
distributions may not exactly coincide there is a strong
dependence on one another). Income, in fact, is only one part of
wealth, the other part being material wealth. Material products
have no conservation law and their monetary value is not constant
(they can be manufactured, destroyed, consumed, etc.). The
distinction between the two will not influence the conclusion
drawn here from the dynamical model of economic interactions we
consider.

Following this line of thought,  we consider here a very simple
model of an open market economy involving both assets exchanges
between individuals and speculative trading. In our non-stationary
economy the total wealth is not conserved due to a random dynamics
which describes the spontaneous growth or decrease of wealth due
to investments in the stock market, housing, ect. It is important
to note that this mechanism corresponds to the effects of an open
market economy where the investments cause the total economy to
growth (more precisely the rich would get richer and the poor
would get poorer). The conservative exchanges dynamic between
individuals redistributes the wealth among people.

We shall show that this kinetic model, introduced in \cite{P},
gives in a suitable asymptotic limit (hereafter called continuous
trading limit) a partial differential equation of Fokker-Planck
type for the distribution of money among individuals.

The equilibrium state of the Fokker-Planck equation can be
computed explicitly and is of Pareto type, namely it is
characterized by a power-law tail for the richest individuals. Our
result can be related to the work \cite{BM}, where a similar
equation has been derived by a mean-field limit procedure.

The mathematical methods we use are close to those used in the
context of kinetic theory for granular flows,  where the limit
procedure is known as quasi-elastic asymptotics  \cite{to1}. We
mention here a similar asymptotic analysis performed on a kinetic
model of a simple market economy with a constant growth mechanism
\cite{Sl}.

The paper is organized as follows. In the next section, we introduce the binary interaction
between agents, which is at the basis of the kinetic model. The main properties of the model are
discussed in section 3. These properties justify the continuous trading limit procedure,
performed in section 4. The limit is illustrated by several numerical results in section 5.

%%%%%%%%%%%%%%%%%%%%%%%%%%%%%%%%%%%%%%%%%%%%%%%%%%%%%%%%%%%%%%%%%%%%%%%%%%

\section{A kinetic model of money asset exchanges}

The goal of a kinetic model of simple market economies, is to describe the evolution of the
distribution of money by means of {\em microscopic} interactions among agents or individuals
which perform exchange of money. Each trade can indeed be interpreted as an interaction where a
fraction of the money changes hands.  We will assume that this wealth after the interaction is
non negative, which corresponds to impose that no debts are allowed. This rule emphasizes the
difference between economic interactions, where not all outcomes are permitted, and the classical
interactions between molecules.

>From a  microscopic view point, the binary interaction is
described by the following rules (we refer to \cite{P} for more
details)
 \bqa \nonumber
w'  & =  & (1 - \gamma) w + \gamma w_* + \eta w\\
\label{trade_rule}
\\[-.25cm]
\nonumber w_*'  & =&   (1 - \gamma) w_* + \gamma w + \eta_* w_*
 \eqa
 where $(w,w_*)$
denote the (positive) money of two arbitrary individuals before
the trade and $(w',w_*')$ the money after the trade. In
(\ref{trade_rule}) we will not allow agents to have debts, and
thus the interaction takes place only if $w' \ge 0$ and $w'_* \ge
0$. In (\ref{trade_rule}) the transaction coefficient $\gamma\in
(0,1)$ is a given constant, while  $\eta$ and $\eta_*$ are random
variables with the same distribution (for example normal) with
variance $\sigma^2$ and zero mean.

 Let us describe the three terms in the right hand side. The first term is related to the marginal
 saving propensity of the agents, the second corresponds to the money transaction, and the last contains
  the effects of
an open economy describing the market returns. Note that since
debts are not allowed,  the total amount of money in the system is
increasing.

This binary interaction model can be also related to some recent
Lotka-Volterra type models \cite{BM}. A similar trade rule where
$\eta$ and $\eta_*$ have the same positive constant value has been
considered recently in \cite{Sl}. In a closed economical system it
is assumed that the total amount of money is conserved ($\eta,
\eta_*\equiv 0$). This conservation law is reminiscent of
analogous conservations which take place in kinetic theory. In
such a situation equations (\ref{trade_rule}) correspond to a
granular gas like model (or to a traffic flow model) where the
stationary state is a Dirac delta centered in the average wealth
(usually referred to as synchronized traffic state in traffic flow
modelling). Thus all agents will end up in the market with exactly
the same amount of money.

 Let  $f(w,t)$ denote the distribution of money $w \in \real_+$ at time $t \ge
0$. By standard methods of kinetic theory, it is easy to show that the time evolution of $f$ is driven by the
following integro-differential equation of Boltzmann type,
 \be \frac{\partial f}{\partial t} =
\int_{\real^2}\int_0^\infty (\beta_{('w , 'w_*) \to (w , w_*)} J f('w)
 f('w_*) - \beta_{(w , w_*) \to (w', w'_*)} f(w) f(w_*) d
 w_*\,d\eta\,d\eta_* \label{eq:boltz}
  \ee
  where $('w,'w_*)$ are the
pre-trading money that generates the couple $(w,w_*)$ after the interaction. In (\ref{eq:boltz}) $J$ is jacobian of
the transformation of $(w,w_*)$ into $(w',w'_*)$ and the kernel  $\beta$ is related to the details of the binary
interaction.

We shall restrict here to a transition rate of the form
\[
\beta_{(w,w_*) \to (w',w_*')}=\mu(\eta)\mu(\eta_*)\Psi(w'\geq 0)\Psi(w'_*\geq 0),
\]
where $\Psi(A)$ is the indicator function of the set $A$, and $\mu(\cdot)$ is a symmetric
probability density with zero mean and  variance $\sigma^2$. The rate function $\beta_{(w',w_*')
\to (w,w_*)}$ characterizes the effects of the open economy through the distribution of the
random variables $\eta$ and $\eta_*$ and takes into account the hypothesis that no-debts are
allowed. The above equation can be included in a more general settings where the trade rule has a
more complex structure including, for example, risk, taxes and subsidies as described in
\cite{P}.

We remark that, for general probability density $\mu(\cdot)$, the
rate function $\beta$ depends on the wealth variables $(w,w_*)$
through the indicator functions $\Psi$. This is analogous to what
happens in the classical Boltzmann equation \cite{Cer}, where the
rate function depends on the relative velocity.  A simplified
situation occurs when the random variables take values on the set
$(-(1-\gamma), 1-\gamma)$. In this case in fact, both $w'\ge 0$
and $w'_* \ge 0$, and the kernel $\beta$ does not depend on the
wealth variables $(w,w_*)$. In this case the kinetic equation
(\ref{eq:boltz}) is the corresponding of the classical Boltzmann
equation for Maxwell molecules \cite{Bob}, which presents several
mathematical simplifications. In all cases however, methods
borrowed from kinetic theory of rarefied gas can be used to study
the evolution of the function $f$.

%%%%%%%%%%%%%%%%%%%%%%%%%%%%%%%%%%%%%%%%%%%%%%%%%%%%%%%%%%%%%%%%%%%%%%%%%%%%%%%%%%%%%%%%%%%%%%%%%%%%%%%%%%%%%%%%%%%%%%%%%

\section{Main properties of the kinetic equation}
We will start our  analysis  by introducing some notations and by
discussing the main properties of the kinetic equation
(\ref{eq:boltz}). Let $Q(f,f)$ denote the interaction integral,
 \be
 Q(f,f)(w) = \int_{\real^2}\int_0^\infty (\beta_{('w , 'w_*) \to (w , w_*)} J f('w)
 f('w_*) - \beta_{(w , w_*) \to (w', w'_*)} f(w) f(w_*) d
 w_*\,d\eta\,d\eta_* .
  \ee

 Let $\bigM_0$ the space of all probability measures in $\real_+$ and by \bq\label{misure} \bigM_{p} =\left\{
\mu \in\bigM_0: \int_{\real_+} |v|^{p}\mu(dv) < +\infty, p\ge 0\right\}, \eq the space of all Borel probability
measures of finite momentum of order $p$, equipped with the topology of the weak convergence of the measures.

Let $\bigF_s(\real_+)$,  be the class of all real  functions on $\real_+$ such that $g(0)= g'(0) =0$, and
$g^{(m)}(v)$ is H\"older continuous of order $\delta$,
 \bq\label{lip} \|g^{(m)}\|_\delta= \sup_{v\not= w} \frac{|g^{(m)}(v) -g^{(m)}(w)|}{
|v-w|^\delta} <\infty,
 \eq
 the integer $m$ and the number $0 <\delta \le 1$ are such that $m+\delta =s$, and
$g^{(m)}$ denotes the $m$-th derivative of $g$.

 In the rest of the paper we will assume that the symmetric probability density $\mu(\eta)$
 which characterizes the transition rate  belongs to $\bigM_{2+\alpha}$, for some $\alpha
 >0$.  Moreover, to simplify the computations, we assume that this density is obtained
 from a given random variable $Y$ with zero mean and unit
 variance, that belongs to $\bigM_{2+\alpha}$. Thus, $\mu$ of variance $\sigma^2$ is the density of
 ${\sigma}Y$. By this assumption, we can easily obtain the dependence on $\sigma$ of the moments of $\mu$.
 In fact, for any $p >2$ such that the $p$-th moment of $Y$ exists,
 \[
\int_{\real}|\eta|^{p}\mu(\eta) d\eta = E\left( \left|{\sigma}Y\right|^{p}\right) =
\sigma^{p}E\left(\left|Y\right|^{p}\right).
\]

By a weak solution of the initial value problem for equation
\rf{eq:boltz}, corresponding to the initial probability density
$f_0(w) \in \bigM_{p}, p >1$ we shall mean any probability density
$f \in C^1(\real_+, \bigM_{p})$ satisfying the weak form of the
equation
 \begin{eqnarray}\label{weak boltz}
&&\frac d{dt}\int_{\real_+} \phi(w)f(w,t)\,dw  = (Q(f,f),\phi) = \nonumber \\
 &&\int_{\real^2} \int_{\real^2_+} \beta_{(w,w_*) \to (w',w'_*)} f(w)
f(w_*) ( \phi(w')-\phi(w)) d w_* d w d \eta\,d\eta_*,
\end{eqnarray}
for $t>0$ and all $\phi \in \bigF_{p}(\real_+)$, and such that for all $\phi \in
\bigF_{p}(\real_+)$
 \bq\label{ic} \lim_{t\to 0} \int_{\real_+} \phi(w)f(w,t)\, dw = \int_{\real_+} \phi(w)f_0(w)\, dw.
  \eq
We can alternatively use the symmetric form
 \begin{eqnarray} \label{weak boltz2}\nonumber
\frac{d}{dt} \int_0^\infty f(w) \phi(w)\,dw &=& \frac12
\int_{\real^2} \int_{\real^2_+}
\beta_{(w,w_*) \to (w',w'_*)}  f(w) f(w_*) \\
\\[-.25cm]
\nonumber &&( \phi(w')+\phi(w'_*)-\phi(w)-\phi(w_*)) dw_*\, dw\,
d\eta\,d\eta_*.
  \end{eqnarray}
Existence of a weak solution to the initial value problem for equation (\ref{eq:boltz}) can be
easily obtained by using methods first applied to the Boltzmann equation \cite{Cer}. The main
difference here is that the classical Boltzmann $H$-theorem, which prevents formation of
concentration, does not hold.

>From (\ref{weak boltz}) (or equivalently from (\ref{weak boltz2}))
conservation of the total number of agents is obtained for
$\phi(w)=1$ which represents the only conservation property
satisfied by the system. The choice  $\phi(w)= w$ is of particular
interest since it gives the time evolution of the total wealth. We
have
 \be
\frac d{dt}\int_{\real_+} w f(w,t)\,dw  =
 \int_{\real^2} \int_{\real^2_+} \beta_{(w,w_*) \to (w',w'_*)} f(w)
f(w_*) (\g (w_* -w) + \eta w ) d w_* d w d \eta\,d\eta_* = \nonumber
 \ee
 \be\label{grow}
-\int_{\real^2} \int_{\real^2_+}\eta \mu(\eta)\mu(\eta_*)\left( 1- \Psi(w'\geq 0)\Psi(w'_*\geq 0)
\right) w f(w) f(w_*)   d w_* d w d \eta\,d\eta_*.
 \ee
 Since $\eta <0$ on the set where $\Psi(w'\geq 0) =0$,
 equation \rf{grow} shows that the total wealth is increasing, unless $w' \ge 0$ and $w'_* \ge 0$
 for all exchanges (in such case the total wealth would be conserved). An upper bound for the total growth
 can be derived considering that $1-\gamma
 +\eta \ge 0$ implies $w'\geq 0$ (respectively $1-\gamma +\eta_* \ge 0$ implies $w'_*\geq 0$).
 Thus,
  \be\label{meglio}
 1- \Psi(w'\geq 0)\Psi(w'_*\geq 0) \le \Psi(1-\gamma +\eta < 0) + \Psi(1-\gamma +\eta_* < 0).
  \ee
Substituting this inequality into \rf{grow} we obtain
 \be
\frac d{dt}\int_{\real_+} w f(w,t)\,dw  \le -\int_{\real^2}\eta \mu(\eta)\mu(\eta_*)\left(
\Psi(1-\gamma +\eta < 0) \right. \nonumber
 \ee
  \be
+\left. \Psi(1-\gamma +\eta_* < 0) \right)d \eta\,d\eta_*\int_{\real_+} w f(w)  d w .
  \ee
By Markov inequality,
 \[
\int_{\real}\mu(\eta_*)d\eta_* \Psi(1-\gamma +\eta_* < 0) \le
\frac{\sigma^{2+\alpha}}{(1-\gamma)^{2+\alpha}},
 \]
while
 \[
 \int_{\real}|\eta| \mu(\eta)d\eta \Psi(1-\gamma +\eta < 0) \le \frac 1{(1-\gamma)^{1+\alpha}}
 \int_{\real}|\eta|^{2+\alpha} \mu(\eta)d\eta \Psi(1-\gamma +\eta < 0) \le
\frac{\sigma^{2+\alpha}}{(1-\gamma)^{1+\alpha }}.
 \]
 Finally,
  \be
\frac d{dt}\int_{\real_+} w f(w,t)\,dw  \le  \frac{\sigma^{2+\alpha}(2-\gamma)}{(1-\gamma)^{2+\alpha}}
\int_{\real_+} w f(w,t)\,dw.
 \ee
We remark that, if the kernel $\beta$ does not depend on the wealth variables, the kinetic model
preserves the total wealth.

Similar bounds for  moments of  order higher than $2$ can be obtained in a similar way. To this
purpose, let us fix $\phi(w) = w^p$ for some $p>2$. Using the same trick as in \rf{grow} we
obtain
 \be
\frac d{dt}\int_{\real_+} w^p f(w,t)\,dw  =
 \int_{\real^2} \int_{\real^2_+} \beta_{(w,w_*) \to (w',w'_*)} f(w)
f(w_*) (|w'|^p -w^p ) d w_* d w d \eta\,d\eta_* = \nonumber
  \ee
  \be
  \int_{\real^2} \int_{\real^2_+} \mu(\eta)\mu(\eta_*) f(w)
f(w_*) (|w'|^p -w^p ) d w_* d w d \eta\,d\eta_* + \nonumber
 \ee
 \be
-\int_{\real^2} \int_{\real^2_+}\mu(\eta)\mu(\eta_*)\left( 1- \Psi(w'\geq 0)\Psi(w'_*\geq 0)
\right)(|w'|^p -w^p ) f(w) f(w_*)   d w_* d w d \eta\,d\eta_* \le \nonumber
 \ee
  \be
  \int_{\real^2} \int_{\real^2_+} \mu(\eta)\mu(\eta_*) f(w)
f(w_*) (|w'|^p -w^p ) d w_* d w d \eta\,d\eta_* + \nonumber
 \ee
 \be\label{growp}
+\int_{\real^2} \int_{\real^2_+}\mu(\eta)\mu(\eta_*)\left( 1- \Psi(w'\geq 0)\Psi(w'_*\geq 0)
\right)w^p  f(w) f(w_*)   d w_* d w d \eta\,d\eta_* .
 \ee
Using again Markov inequality we obtain
 \be
 \int_{\real^2} \int_{\real^2_+}\mu(\eta)\mu(\eta_*)\left( 1- \Psi(w'\geq 0)\Psi(w'_*\geq 0)
\right)w^p  f(w) f(w_*)   d w_* d w d \eta\,d\eta_* \le \nonumber
 \ee
 \be
\int_{\real^2} \mu(\eta)\mu(\eta_*)\left( \Psi(1-\gamma +\eta < 0)+ \Psi(1-\gamma +\eta_* < 0)
\right)d \eta\,d\eta_*\int_{\real_+} w^p f(w)  d w \le \nonumber
  \ee
  \be
 2 \frac{\sigma^{2+\alpha}}{(1-\gamma)^{2+\alpha}}\int_{\real_+}w^p  f(w)  d w.
 \ee

Moreover, we can write
 \[
 |w'|^p = w^p + pw^{p-1}(w' - w) + \frac 12 p(p-1)|\tilde w|^{p-2}(w'-w)^2,
 \]
 where, for some $0 \le \theta \le 1$
 \[
 \tilde w = \theta w' +(1 -\theta)w .
  \]
 Hence,
 \be
  \int_{\real^2} \int_{\real^2_+} \mu(\eta)\mu(\eta_*) f(w)
f(w_*) (|w'|^p -w^p ) d w_* d w d \eta\,d\eta_* = \nonumber
 \ee
 \be
  \int_{\real^2} \int_{\real^2_+} \mu(\eta)\mu(\eta_*) f(w)
f(w_*) (pw^{p-1}(w' - w) + \frac 12 p(p-1)|\tilde w|^{p-2}(w'-w)^2 ) d w_* d w d \eta\,d\eta_*
\le \nonumber
 \ee
 \be
 \frac 12 p(p-1) \int_{\real^2} \int_{\real^2_+} \mu(\eta)\mu(\eta_*) f(w)
f(w_*)  |\tilde w|^{p-2}(w'-w)^2  d w_* d w d \eta\,d\eta_* .
 \ee
 In fact,
\be
  \int_{\real^2} \int_{\real^2_+} \mu(\eta)\mu(\eta_*) f(w)
f(w_*)w^{p-1}(w' - w) d w_* d w d \eta\,d\eta_* = \nonumber
 \ee
\be
  \int_{\real^2} \int_{\real^2_+} \mu(\eta)\mu(\eta_*) f(w)
f(w_*)w^{p-1}(\gamma(w_* - w)+\eta w) d w_* d w d \eta\,d\eta_* = \nonumber
 \ee
\be
 \gamma \int_{\real^2_+}  f(w)
f(w_*)w^{p-1}(w_* - w) d w_* d w  \le 0.
 \ee
The last bound follows by H\"older inequality. In fact,
 \[
\int_{\real^2_+} f(w) f(w_*)w^{p-1}w_*  d w_* d w =\int_{\real_+} f(w)w^{p-1}d w \int_{\real_+}f(w_*)w_*  d w_* \le
 \]
 \[
\left(\int_{\real_+} f(w)w^{p}d w\right)^{(p-1)/p}\left( \int_{\real_+}f(w_*)w_*^p  d w_*\right)^{1/p} =
\int_{\real_+} f(w)w^{p}d w.
 \]
Finally, since
 \[
|\tilde w|^{p-2} = |(1 -\theta \gamma)w + \theta \gamma w_* + \theta\eta w|^{p-2} \le c_p\left[
w^{p-2} +  \gamma^{p-2} w_*^{p-2} + |\eta |^{p-2}w^{p-2}\right],
 \]
 we obtain
 \be
  \int_{\real^2} \int_{\real^2_+} \mu(\eta)\mu(\eta_*) f(w) f(w_*)  |\tilde w|^{p-2}(w'-w)^2 d
w_* d w d \eta\,d\eta_* \le \nonumber
 \ee
 \be
  c_p\int_{\real^2} \int_{\real^2_+} \mu(\eta)\mu(\eta_*) f(w) f(w_*)\left[
w^{p-2} +  \gamma^{p-2} w_*^{p-2} + |\eta |^{p-2}w^{p-2}\right]\cdot \nonumber
 \ee
  \be
 \cdot\left[ \gamma^2(w_* -w)^2 + \eta^2 w^2 +2\eta w(w_*-w)\right]
   dw_* d w d \eta\,d\eta_* \le \nonumber
 \ee
  \be\label{grow2}
 c_pA_p(\sigma,\gamma)\int_{\real_+}w^p  f(w)  d w,
 \ee
 where
  \be
 A_p(\sigma,\gamma)= \sigma^2( 1 + \gamma^{p-2}) +\sigma^p +2\gamma^2
 +2 \gamma^p +2\gamma^2\sigma^{p-2}.
  \ee
  Grouping all estimates together, we finally obtain the bound for
  the moments
   \be\label{grow3}
\frac d{dt}\int_{\real_+} w^p f(w,t)\,dw  \le \left( \frac 12 p(p-1)c_pA_p(\sigma,\gamma)+ 2
\frac{\sigma^{2+\alpha}}{(1-\gamma)^{2+\alpha}}\right) \int_{\real_+}w^p  f(w) d w.
  \ee

%%%%%%%%%%%%%%%%%%%%%%%%%%%%%%%%%%%%%%%%%%%%%%%%%%%%%%%%%%%%%%%%%%%%%%%%%%

\section{The continuous trading limit}

The previous analysis shows that in general it is quite difficult
both to study in details the evolution of the wealth function, and
to describe the asymptotic behavior. As is usual in kinetic
theory, particular asymptotics of the equation result in
simplified models (generally of Fokker-Planck type), for which it
is relatively easier to find steady states, and to prove their
stability. In order to give a physical basis on these asymptotics,
it is relevant to discuss in some detail the interaction rule
(\ref{trade_rule}). To skip inessential difficulties, that as we
will see later on do not change the structure of the limit
equation, we suppose here that the random variables take values on
the set $(-(1-\gamma), 1-\gamma)$. As remarked in the previous
section, in this case  both $w'\ge 0$ and $w'_* \ge 0$, and the
kernel $\beta$ does not depend on the wealth variables $(w,w_*)$.
Let us denote by $E(X)$ the mathematical expectation of the random
variable $X$. Then the following properties follow from
(\ref{trade_rule})
 \be\label{mean}
 E[w' + w'_*] = w+w_*, \quad   \quad E[w' - w'_*] = (1-2\gamma)(w-w_*).
 \ee
 The first equality in \rf{mean} describes the property of mean conservation of wealth. The
 second refers to the tendency of the trade to decrease (in mean) the distance between wealths after the
interaction. This tendency is a universal consequence of the rule
(\ref{trade_rule}), in that it holds whatever distribution one
assigns to $\mu$, namely to the random variable which accounts for
the effects of the market returns in an open economy. This
universality is false for the first equality, which in general has
to be substituted by the inequality
  \be\label{mean1}
 E[w' + w'_*] \ge  w+w_* .
 \ee
 Hence, in general the effects of the market returns in an open economy account for an increasing
 of the microscopic mean wealth.

 The second property in \rf{mean} is analogous to the similar one that holds in a collision between molecules in a
 granular gas. There the quantity $e=2\gamma$ is called "coefficient of restitution", and
 describes the peculiar fact that energy is dissipated \cite{to1}. If we want to consider the situation in which most of the trades
corresponds to a very small exchange of money ($\gamma \to 0$), and at the same time we want to
maintain both properties \rf{mean} at a macroscopic level, we have to pretend that
 \be\label{11}
\int_{\real_+^2}(w+w_*)f(w)f(w_*) dw dw_* = 2\int_{\real_+} w f(w) dw = 2m(t)
 \ee
 remains constant, while
 \be\label{22}
\int_{\real_+^2}(w-w_*)^2f(w)f(w_*) dw dw_* = A_f(t)
 \ee
 is varying with time, and decays to zero when the market returns are not present (i.e. $\sigma=0$).

When the kernel $\beta$ does not depend on the wealth variables,
\rf{grow} implies that $m(t)=m_0$. Moreover, explicit computations
show that $A_f(t)$ varies with law
 \be\label{33}
 \frac {dA_f(t)}{dt} = -4\left[ \gamma(1-\gamma)- \frac{\sigma^2}2 \right]A_f(t) + 2\sigma^2 m^2.
 \ee
 Hence, when $\sigma = 0$, $A_f(t)$ decays exponentially to zero.
It is now evident that \rf{11} is satisfied for any value of
$\gamma$ and $\sigma$, while \rf{22} looses its meaning
as $\gamma$ an $\sigma$ tend to zero. Of course, one can rewrite
\rf{33} as \be\label{44}
 \frac {dA_f(t)}{dt} = -4\gamma\left[ \left( 1-\gamma- \frac{\sigma^2}{2\gamma}\right) A_f(t)
 -
 \frac{\sigma^2}{2\gamma} m^2\right].
 \ee
Hence, if we set
 \be\label{resc}
 \tau = \gamma t, \quad g(w,\tau) = f(w,t),
 \ee
 which implies $f_0(w) = g_0(w)$, we obtain
 \be
\label{45}
 \frac {dA_g(\tau)}{d\tau} = -4\left( 1-\gamma- \frac{\sigma^2}{2\gamma}\right) A_g(\tau) +
 2\frac{\sigma^2}\gamma m^2.
 \ee
 Letting now both $\gamma \to 0$ and $\sigma \to 0$ in such a way that $\sigma^2/\gamma =
 \lambda$, \rf{45} becomes in the limit
 \be
\label{46}
 \frac {dA_g(\tau)}{d\tau} = -\left( 4- 2\lambda\right) A_g(\tau) +
 2\lambda m^2.
 \ee
This formal argument shows that the value of the ratio
$\sigma^2/\gamma$ is of paramount importance to get asymptotics
which maintain memory of the microscopic interactions. In
particular, while for $\lambda <2$  $A_g(\tau)$ converges to the
finite value $ \bar A_g = {\lambda m^2}/(2-\lambda)$ , $A_g(\tau)$
diverges to infinity as time goes to infinity when $\lambda \ge
2$.

 In the remaining
of this section, we shall present a rigorous derivation of a Fokker-Planck model from the
Boltzmann equation for the wealth function $g(w, \tau)$, when both $\gamma \to 0$ and $\sigma \to
0$ in such a way that $\sigma^2/\gamma \to
 \lambda $. This derivation, which is similar to the quasi-elastic limit of
granular gases, is of major relevance for the study of the asymptotic equilibrium states of the kinetic model.
First, we show how the Fokker-Planck equation comes out for the simpler case of a kernel which does not depend on
the wealth variables. Second, we extend the result to a general kernel.

The scaled density  $g(v,\tau)=f(v,t)$ satisfies the weak form
 \be\label{evol}
 \frac{d}{d\tau} \int_0^\infty g \phi\,dw
= \frac1{\gamma} \int_{\real^2} \int_{\real^2_+} \mu(\eta)\mu(\eta_*) g(w) g(w_*) (
\phi(w')-\phi(w)) d w_* d w d \eta\,d\eta_*.
 \ee
Given $0 <\delta \le \alpha$, let us set $\phi \in \bigF_{2+\delta}(\real_+)$.

By (\ref{trade_rule}),
$$
w'   - w =   \gamma (w_* - w) + \eta w .
$$
Then, if we use a second order Taylor expansion of $\phi$ around $w$
 $$
\phi(w') - \phi(w) = (\gamma (w_* - w) + \eta w) \phi'(w) + {1 \over 2} (\gamma (w_* - w) + \eta
w)^2 \phi''(\tilde w),
 $$
where, for some $0 \le \theta \le 1$
 \[
 \tilde w = \theta w' +(1 -\theta)w .
  \]

 Inserting this expansion in the collision operator,  we get
 \be
\frac{d}{d\tau} \int_0^\infty g \phi\,dw = \frac1{\gamma}\int_{\real^2}\int_{\real^2_+}\mu(\eta)\mu(\eta_*)
 [  (\gamma (w_* - w) + \eta w) \phi'(w)\nonumber +
  \ee
  \be\label{form}
  + {1 \over 2}
(\gamma (w_* - w) + \eta w)^2 \phi''(w) ]
 g(w_*)g(w)  d w_*\,d w\,d \eta\,d\eta_*  + R(\gamma, \sigma),
 \ee
where
 \begin{eqnarray*}
 R(\gamma, \sigma) &=
&\frac 1{2\gamma}\int_{\real^2}\int_{\real^2_+}\mu(\eta)\mu(\eta_*) (\gamma (w_* - w) + \eta
w)^2\cdot \\&& \cdot \left( \phi''(\tilde w)- \phi''( w)\right)
 g(w_*)g(w)  d w_*\,d w\,d \eta\,d\eta_*
 \end{eqnarray*}
Since $\phi \in \bigF_{2+\delta}(\real_+)$, and $|\tilde w - w |= \theta|w'-w|$
 \be\label{rem}
 \left| \phi''(\tilde w)- \phi''( w)\right| \le \| \phi''\|_\delta |\tilde w - w |^\delta \le
 \| \phi''\|_\delta |w' - w |^\delta .
  \ee
  Hence
 \begin{eqnarray*}
|R(\gamma, \sigma)| &\le& \frac{\|
\phi''\|_\delta}{2\gamma}\int_{\real^2}\int_{\real^2_+}\mu(\eta)\mu(\eta_*)\cdot \\&& \cdot
|\gamma (w_* - w) + \eta w|^{2+\delta} g(w_*)g(w)  d w_*\,d w\,d \eta\,d\eta_*
 \end{eqnarray*}
By virtue of the inequality
 \[
 |\gamma (w_* - w) + \eta w|^{2+\delta} \le 4^{1+\delta}
 \left( |\gamma w_*|^{2+\delta} +|\gamma w|^{2+\delta}+| \eta w|^{2+\delta}\right),
 \]
we finally obtain the bound
 \be\label{to0}
|R(\gamma, \sigma)| \le 2^{1+2\delta}{\| \phi''\|_\delta}\left( 2\gamma^{1+\delta} + \frac
1\gamma \int_{\real}|\eta|^{2+\delta}\mu(\eta) d\eta \right)\int_{\real_+}w^{2+\delta}g(w) dw
 \ee
 Since $\mu$ is a probability density with zero mean and
$\lambda\gamma$ variance, and $\mu$ belongs to $\bigM_{2+\alpha}$, for  $\alpha
 >\delta$,
 \[
\int_{\real}|\eta|^{2+\delta}\mu(\eta) d\eta = E\left(
\left|\sqrt{\lambda\gamma}Y\right|^{2+\delta}\right) =
(\lambda\gamma)^{1+\delta/2}E\left(\left|Y\right|^{2+\delta}\right),
\]
and $E\left(\left|Y\right|^{2+\delta}\right)$ is bounded. Using this equality into \rf{to0} one
shows that $R(\gamma, \sigma)$ converges to zero as $\gamma \to 0$, if
$\int_{\real_+}w^{2+\delta}g(w,\tau)$ remains bounded at any fixed time $\tau >0$, provided the
same bound holds at time $\tau =0$. By virtue of bound \rf{grow3}, taking $\sigma^2 =
\lambda\gamma$, we obtain
   \be\label{grow4}
\frac d{d\tau}\int_{\real_+} w^{2+\delta} g(w,\tau)\,dw  \le \nonumber
 \ee
  \be
\frac 1\gamma\left( \frac 12 p(p-1)c_{2+\delta}A_{2+\delta}(\sqrt{\lambda\gamma},\gamma)+ 2
\frac{(\lambda\gamma)^{1+\alpha/2}}{(1-\gamma)^{2+\alpha}}\right) \int_{\real_+}w^p  g(w,\tau) d
w.
  \ee
Since the lower order term in $A_{2+\delta}(\sqrt{\lambda\gamma},\gamma)$ is $\lambda\gamma(1-
\gamma)^\delta$, it follows that the boundedness of the moment holds independently of the value
of $\gamma$. Therefore, at any fixed time $\tau$
 \be\label{to1}
|R(\gamma, \sigma)|(\tau) \le 2^{1+2\delta}{\| \phi''\|_\delta}\left( 2\gamma^{1+\delta} +
\lambda^{2+\delta}\gamma^{\delta/2}E\left(\left|Y\right|^{2+\delta}
 \right)\right)C\tau\int_{\real_+}w^{2+\delta}g_0(w).
dw
 \ee
Hence, the remainder $R(\gamma, \sigma)$ converges to zero as both $\gamma$ and $\sigma$ converge
to zero, in such a way that $\sigma^2 = \lambda\gamma$. Within the same scaling,
 \be
 \lim_{\gamma \to 0}
\frac1{\gamma}\int_{\real^2}\int_{\real^2_+}\mu(\eta)\mu(\eta_*)
 [  (\gamma (w_* - w) + \eta w) \phi'(w) + \nonumber
 \ee
 \be
 + {1 \over 2}
(\gamma (w_* - w) + \eta w)^2 \phi''(w) ]
 g(w_*)g(w)  d w_*\,d w\,d \eta\,d\eta_* = \nonumber
  \ee
  \be\label{FP1}
\int_{\real_+}\left[
 (m - w)  \phi'(w)  +
   \frac\lambda{2} w^2 \phi''(w)\right]
 g(w) d w
  \ee
 The right-hand side of \rf{FP1} is nothing but the weak form of the Fokker-Planck equation
 \be\label{FP}
 \frac{\partial g}{\partial \tau} = \frac \lambda 2\frac{\partial^2 }{\partial w^2}\left(w^2
 g\right) + \frac{\partial }{\partial w}\left((w -m)
 g\right).
 \ee
The limit Fokker-Planck equation can be rewritten as
 {\be \frac{\partial g}{\partial \tau} =
\frac{\partial}{\partial w} [ (w-m) + \frac{\lambda}{2} w ) g + \frac{\lambda}{2} w
\frac{\partial}{\partial w} (w  g) ].
 \ee

The general case of a wealth-variables depending rate function can be easily obtained from the previous
computations. With respect to formula \rf{form} we have to consider two more terms in the remainder. The first one
comes out from the possibility to have negative wealth variables as outcome  of the interaction. This term
reads
 \be\label{rem1}
R_1(\gamma, \sigma) = - \frac 1\gamma \int_{\real^2} \int_{\real^2_+} \mu(\eta)\mu(\eta_*)\cdot\nonumber
 \ee
 \be
\cdot \left( 1- \Psi(w'\geq 0)\Psi(w'_*\geq 0) \right)(\phi(w')-\phi(w)) g(w) g(w_*)   d w_* d w d \eta\,d\eta_*.
 \ee
 The integrand in \fer{rem1} coincides with the integrand in the right-hand side of \fer{evol}, multiplied by the
 factor $1- \Psi(w'\geq 0)\Psi(w'_*\geq 0)$. Using \fer{meglio} shows that we have to integrate the probability
 density functions on a bounded domain, which gives a better behavior as $\gamma, \sigma \to 0$. In more details,
 while
 \[
\int_{\real}\eta^2  \mu(\eta)\, d\eta = \sigma^2,
 \]
 \[
\int_{\{|\eta| > 1-\gamma \}}\eta^2  \mu(\eta)\, d\eta \le \frac 1{(1-\gamma)^\alpha} \int_{\{|\eta| > 1-\gamma
\}}\eta^{2+\alpha} \mu(\eta)\, d\eta \le \sigma^{2+\alpha}E(Y^{2+\alpha}).
 \]
Thus, we can use the same expansion  following formula \fer{evol} to
 conclude that  $R_1(\gamma, \sigma)$ converges to zero as $\gamma \to 0$.

 The second remainder we have to take into account comes out from the fact that in general the total wealth $m(t)$
is not constant in time. Thus we have the additional term
 \be
R_2(\gamma, \sigma) = - \frac 1\gamma \int_{\real^2}  \mu(\eta)\mu(\eta_*)\left( \int_{\real_+}w_*g(w_*,\tau)
\,dw_* - m\right) g(w,\tau)\phi'(w) d w_* d w d \eta\,d\eta_*.
 \ee
On the other hand, thanks to \rf{evol}
 \[
\left| \int_{\real_+}w_*g(w_*,\tau) \,dw_* - m\right| \le \left[
\exp\left\{\frac{\sigma^{2+\alpha}(2-\gamma)}{(1-\gamma)^{2+\alpha}}\tau \right\} - 1 \right] m,
 \]
 which implies the convergence to zero of the remainder as $\gamma \to 0$.
 Hence we proved

\begin{theorem}
 Let the probability density $f_0 \in \bigM_p$, where $p= 2+ \delta$ for some $\delta>0$, and let the symmetric
 random variable $Y$ which characterizes
 the kernel have a density in $\bigM_{2+\alpha}$, with $\alpha > \delta$. Then, as $\gamma \to 0$,
 $\sigma \to 0$ in such a way that $\sigma^2 = \lambda\gamma$  the weak solution to the Boltzmann equation for the scaled density
 $g_\gamma(v,\tau)=f(v,t)$, with $\tau = \gamma t$
converges, up to extraction of a subsequence, to a probability density $g(w,\tau) $. This density  is a weak
solution of the  Fokker-Planck equation \rf{FP}, and it is such that the mean wealth is conserved.
\end{theorem}

The stationary state of the Fokker-Planck equation can be directly computed and, by assuming for simplicity
$$ m=\int_{\real_+}
f(w,t)\,dw=1,$$ it can be written as {
 \be\label{equi}
g_\infty(w)=\frac{(\mu-1)^\mu}{\Gamma(\mu)}\frac{\exp\left(-\frac{\mu-1}{w}\right)}{w^{1+\mu}}
\ee} where {$$ \mu = 1 + \frac{2}{\lambda} >1.
$$}
Therefore the stationary distribution exhibits a Pareto power law
tail for large $w$'s.

Note that this equation is essentially the same Fokker-Planck
equation derived from a Lotka-Volterra interaction in \cite{BM,
So, SMBSR}.

%%%%%%%%%%%%%%%%%%%%%%%%%%%%%%%%%%%%%%%%%%%%%%%%%%%%%%%%%%%%%%%%%%%%%%%%%%

\section{Numerical results}

In this section, we shall present some numerical test and we shall compare the stationary results
obtained by using Monte Carlo simulation of the kinetic model with the stationary state of the
Fokker-Planck model. We start from a situation where all individuals shares the same money. At
each iteration, we choose randomly two individual. Then, the parameter $\eta$ is choosen
accordingly to a normal centered law (one for each individual) and the trade is performed if it
is admissible i.e. when the new money of each individual remain positive. We report the results
for the density { $\tilde{f}(w,t)=m(t)f(m(t)w,t))$} where $m(t)$ is the (variable in time) total
amount of money in the system. Note that the density $\tilde f$ is normalized to have the first
moment equal to one.

We use $N=2000$ individuals and perform several iterations until a
stationary state is reached. The distribution is then averaged
over the next 250 iterations.

We plot the distribution function for different values of
$(\gamma, \sigma^2)=\{(0.1,0.2), (0.01,0.02)\}$.

%\begin{center}
%\begin{tabular}{c|c|c|}
%Figure & $\gamma$ & $\sigma^2$ \\
%\hline
%1 & 0.1 & 0.2 \\
%2 & 0.01 & 0.02 \\
%\end{tabular}
%\end{center}

\begin{figure}[t]
\begin{center}
\includegraphics[scale=.4]{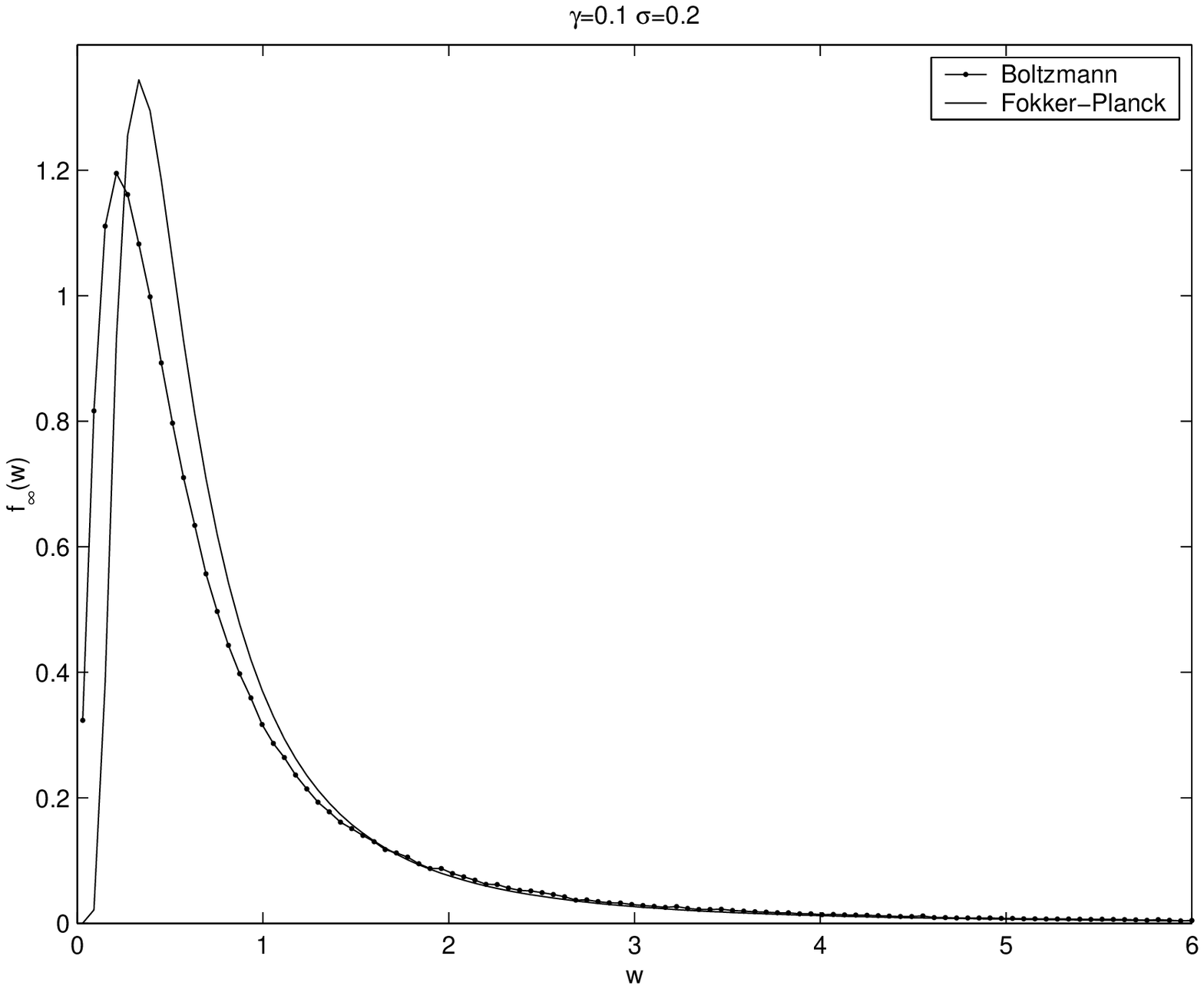}
\includegraphics[scale=.4]{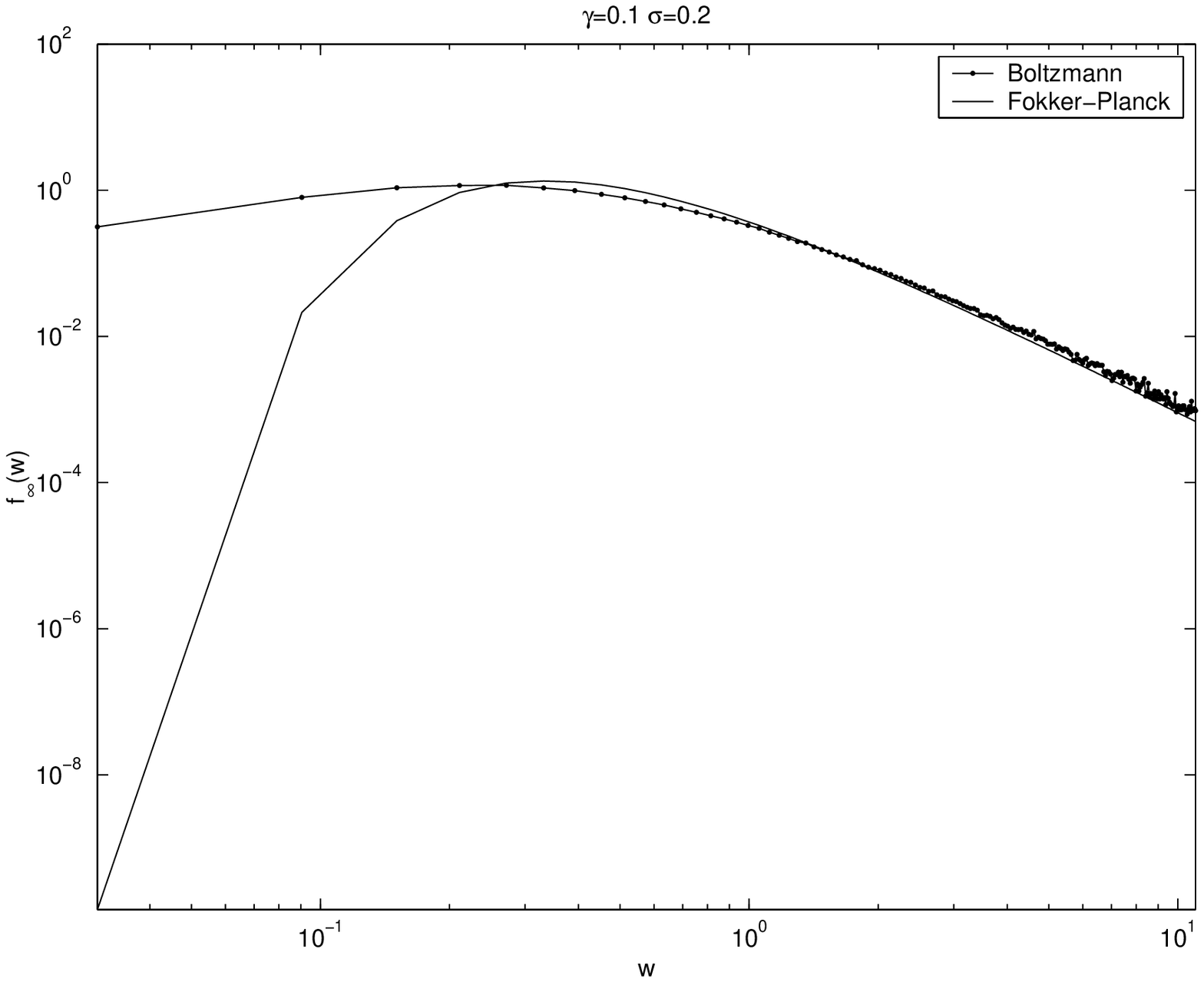}
\end{center}
\caption{Asymptotic behavior of the Fokker-Planck model and the
Boltzmann model for $\mu=2.0$, $\gamma=0.1$ and $\sigma=0.2$.
Figure on the right is in loglog-scale.} \label{fig:1}
\end{figure}

\begin{figure}[t]
\begin{center}
\includegraphics[scale=.4]{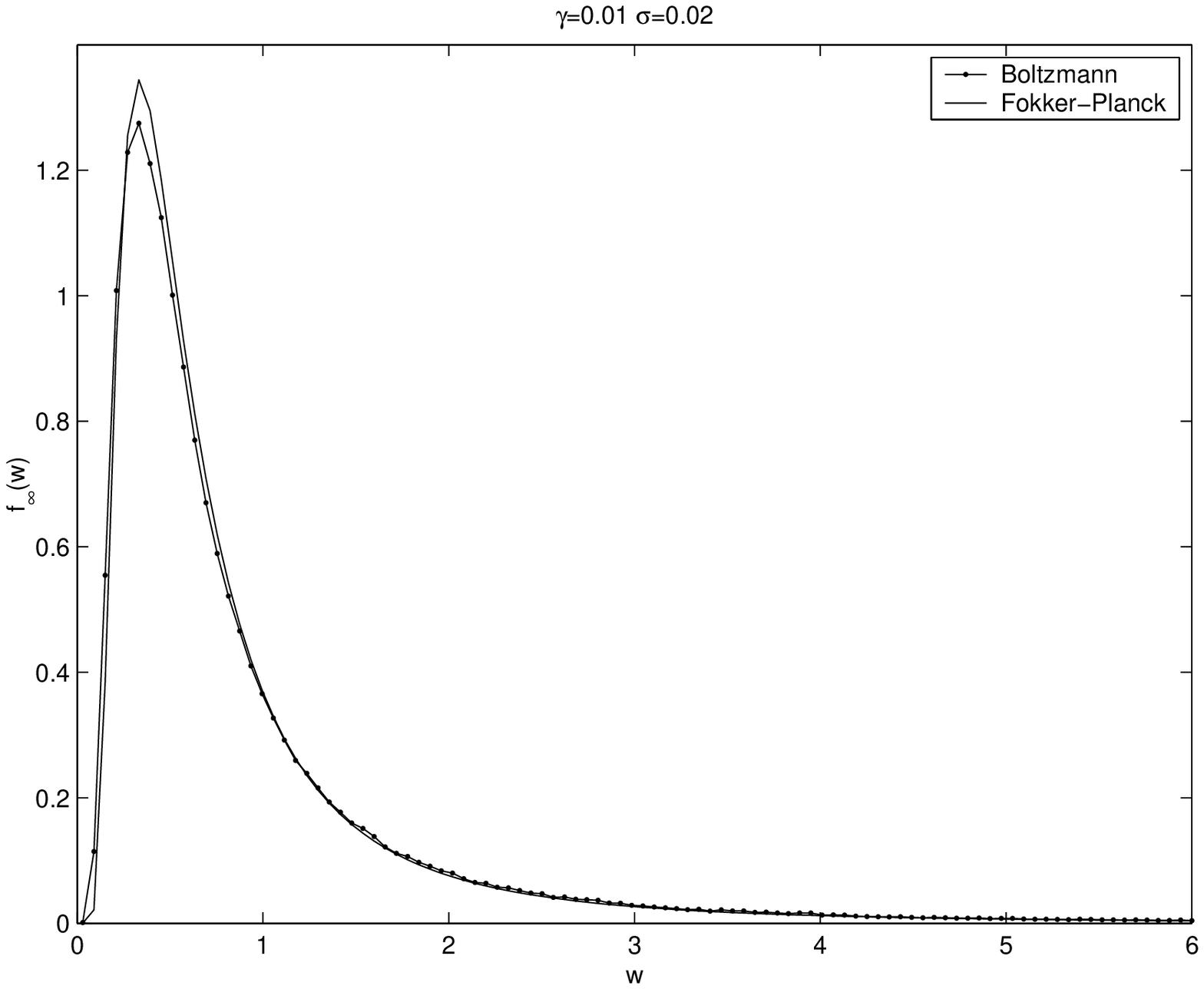}
\includegraphics[scale=.4]{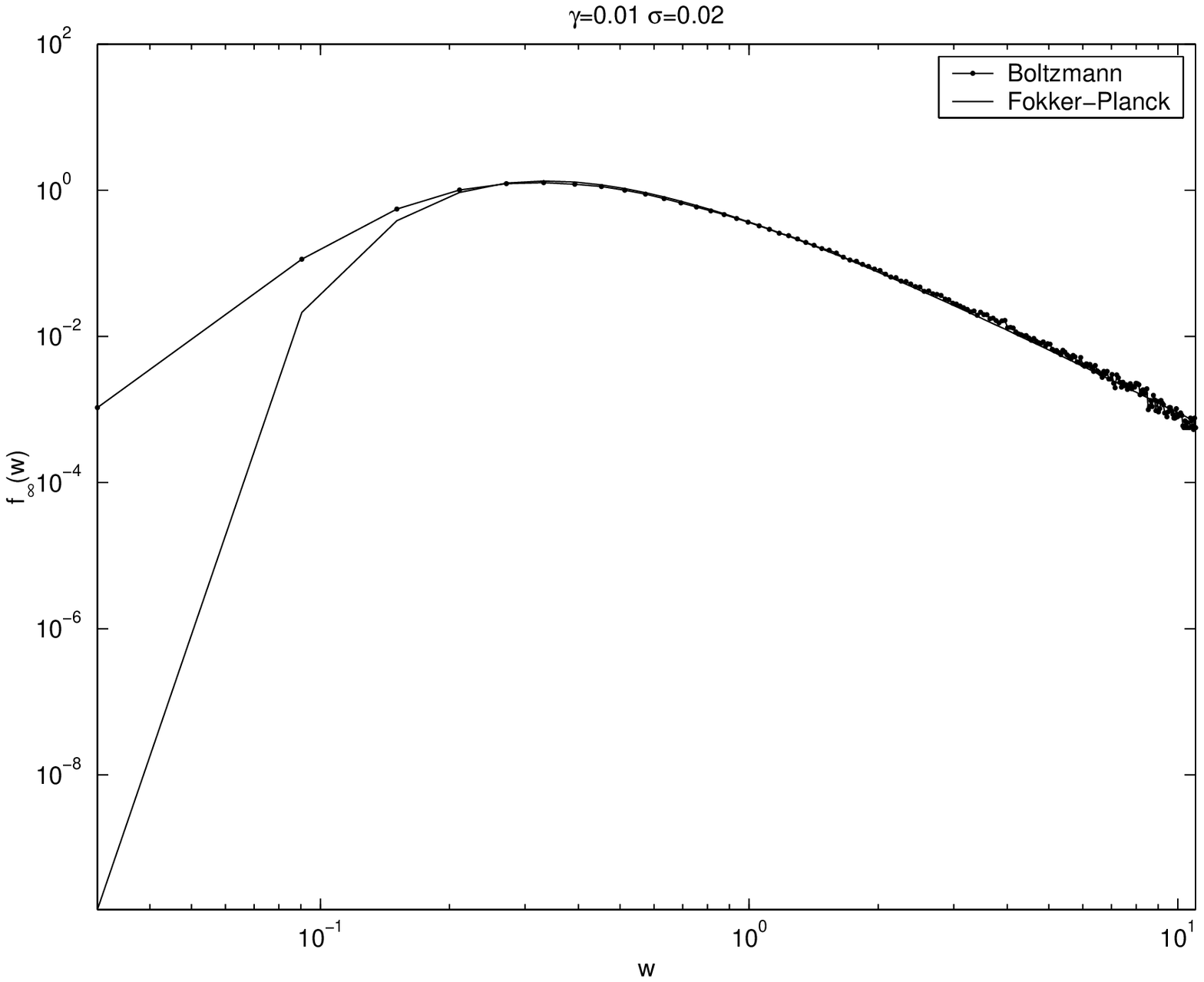}
\end{center}
\caption{Asymptotic behavior of the Fokker-Planck model and the
Boltzmann model for $\mu=2.0$, $\gamma=0.01$ and $\sigma=0.02$.
Figure on the right is in loglog-scale.} \label{fig:2}
\end{figure}

The values for $\gamma$ and $\sigma^2$ are such that we have a
fixed $\lambda=2$ corresponding to a coefficient $\mu=2$ which is
a realistic value for the income distribution of capitalistic
economies \cite{BM, SMBSR}. As prescribed from our theoretical
analysis we observe that the equilibrium distribution converges
toward the function $f_{\infty}$ as both $\gamma$ and $\sigma^2$
go to 0, with $\sigma^2/\gamma = 2$ (see Figures \ref{fig:1} and
\ref{fig:2}).

\begin{figure}[ht]
\begin{center}
\includegraphics[scale=.4]{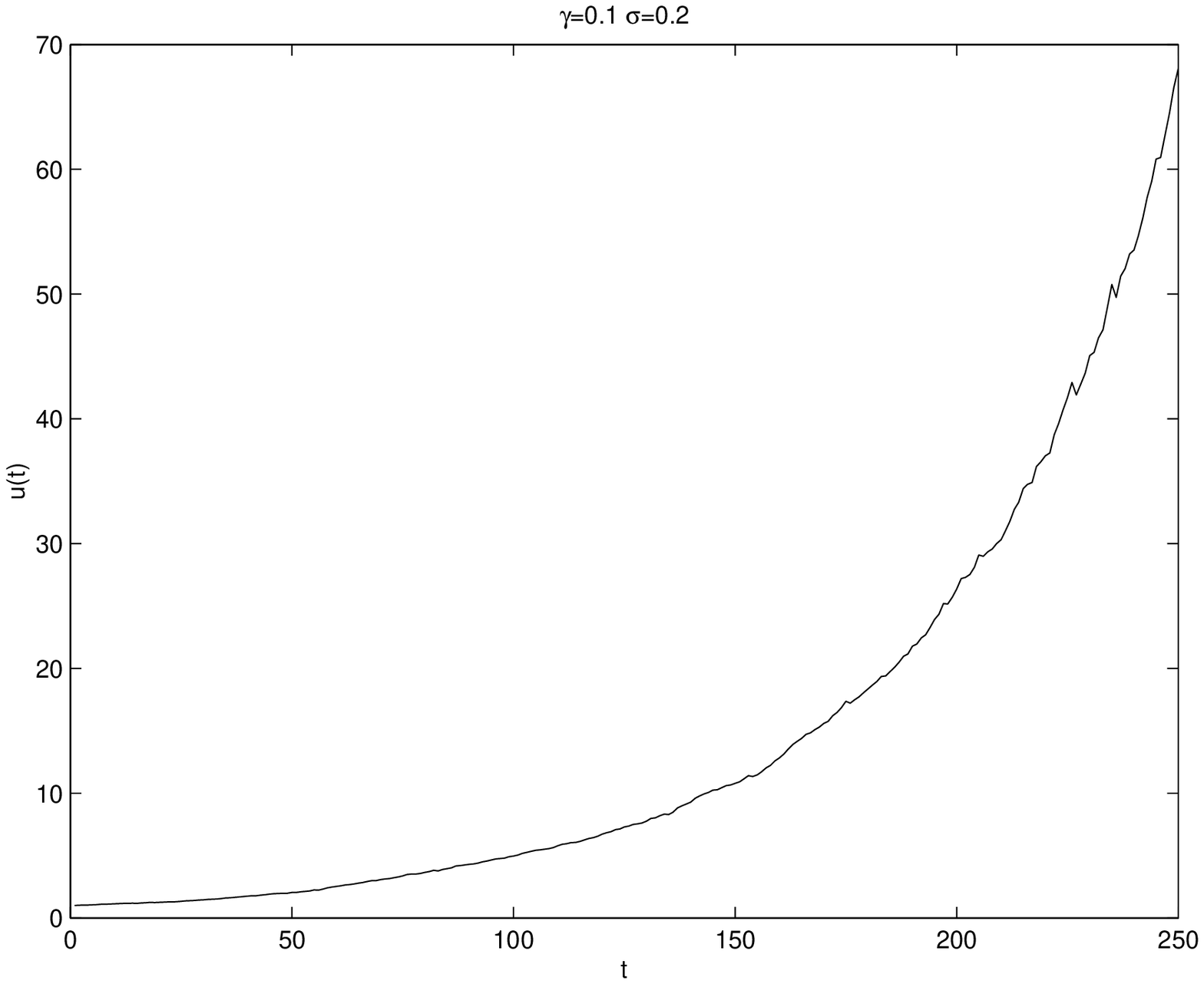}
\includegraphics[scale=.4]{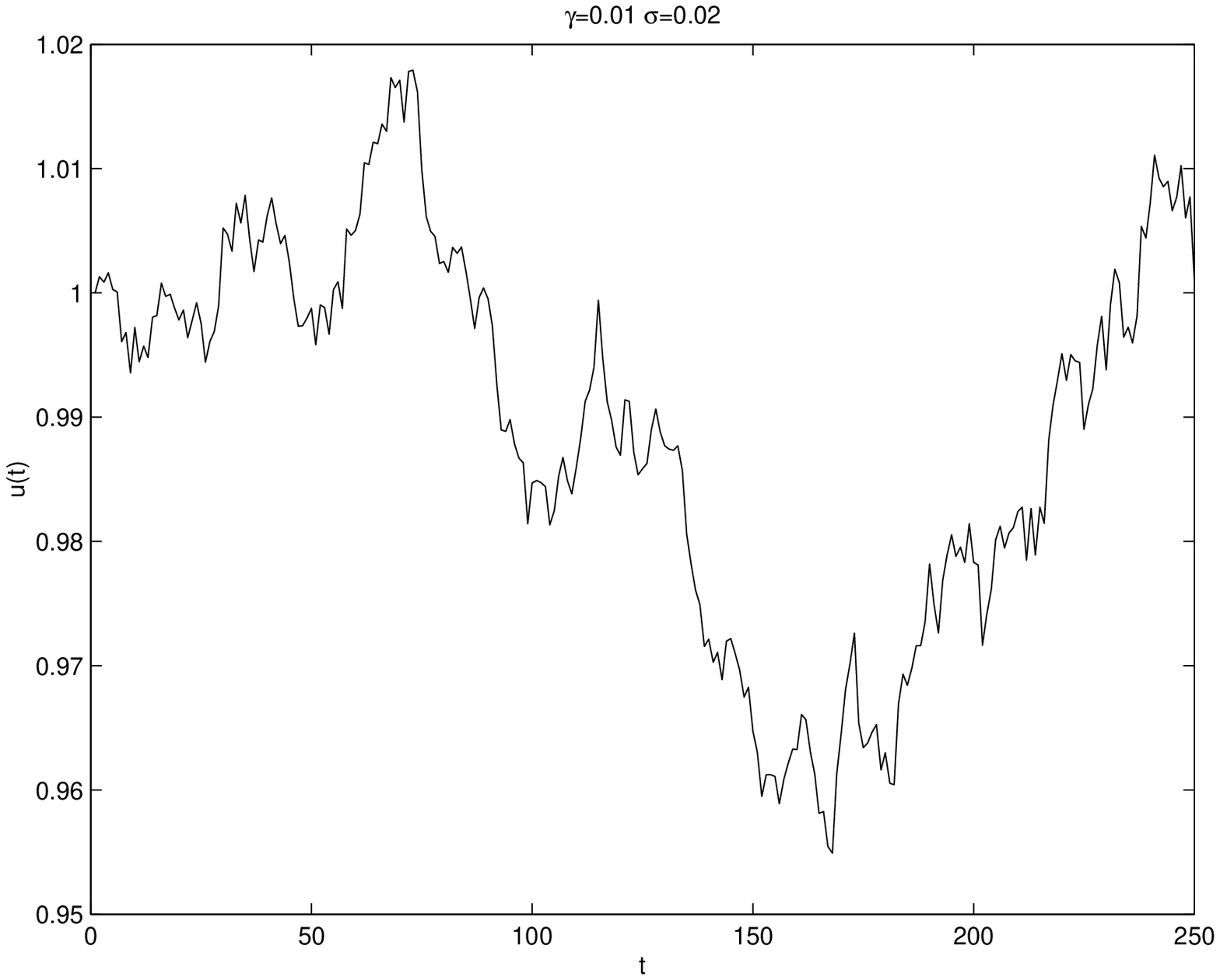}
\end{center}
\caption{Growth of the total amount of money in the market
corresponding to the data of the previous figure.} \label{fig:3}
\end{figure}

The corresponding behaviors of the total amount of money obtained
with $N=10000$ agents are given in Figure \ref{fig:3}. Note that,
as expected, asymptotically the exponential growth vanishes and
the model preserves the total amount of money.

%Let us plot the evolution of the second moment (which can be
%interpreted as a temperature)
%$$
%T = \int_0^\infty f(w) (w-\bar w) dw
%$$
%where $\bar w$ is the constant mean money (equal to one in the
%simulation). It starts from 0 (since all individuals have the same money)
%at the begining and increases quasi-linearly in the begining but it
%is very noisy...

%\includegraphics{temp.pdf}

\section{Conclusions}

We introduced and discussed a nonlinear kinetic model for a simple
market economy, which is based on binary exchanges of money and
speculative trading. We showed that at suitably large times, in
presence of a large number of trades in which agents exchange a
small amount of money, the nonlinear kinetic equation is
well-approximated by a linear Fokker-Planck type equation, which
admits a stationary steady state with Pareto tails. Convergence
towards a similar steady state is shown numerically for the
solution of the kinetic model after a suitable normalization which
guarantees the conservation of the mean wealth. Our analysis
enlightens both the range of validity of the Fokker-Planck
equation \fer{FP}, and why in the continuous trading limit the
mean wealth, which is increasing for the kinetic model, remains
constant.
%The same Fokker-Planck equation was obtained in
%\cite{BM} as the mean field limit of a stochastic dynamical
%equation, as well as in \cite{So, SMBSR} in the context of a
%generalized Lotka-Volterra equation.
%These derivations, however,
%have linear models as starting points.
The present derivation takes advantages of the deep similarities between the trade rule \fer{trade_rule} and a
molecular dissipative collision. However, it is expected that the asymptotic analysis discussed in this paper is
not restricted to the trade rule \fer{trade_rule}, so that different microscopic interactions could be treated as
well. It is clear that the analogy between the trade rule \fer{trade_rule} and a one-dimensional molecular
dissipative collision suggests in a natural way the continuous trading asymptotic which is well-understood in
kinetic theory \cite{to1} with the name of quasi-elastic asymptotic procedure. Furthermore, the formation of
overpopulated energy tails for large times in the kinetic model is in accord with the analogous result valid for
the Boltzmann equation for a dissipative granular Maxwellian gas \cite{BMP, EB1, EB2}. Pushing further these
analogies may help to clarify many aspects of large-time behavior of market economies.

\bigskip \noindent {\bf Acknowledgment:} The authors acknowledge support from the IHP
network HYKE ``Hyperbolic and Kinetic Equations: Asymptotics, Numerics, Applications''
HPRN-CT-2002-00282 funded by the EC., and from the Italian MIUR, project ``Mathematical Problems
of Kinetic Theories''. L. Pareschi thanks the Department of Mathematics and Applications,
Mathematical Physics of Orleans (MAPMO) for the kind hospitality during his visit.

\end{document}